\documentclass[onecolumn,11pt]{article}
\usepackage[top=1in, bottom=1in, left=1in, right=1in]{geometry}
\usepackage{amsmath,amsfonts,amscd,amssymb}
\usepackage{caption}
\usepackage{epstopdf}
\usepackage{overpic}
\usepackage{cancel}
\usepackage{rotating}
\usepackage{url}
\usepackage{color}
\usepackage{rotating}
\usepackage{multirow}
\usepackage{wrapfig}
\usepackage{mathtools}
\usepackage{subeqnarray}
\usepackage{setspace}
\usepackage[numbers]{natbib}

\usepackage{graphicx}
\usepackage{floatrow}
\usepackage{subfig}

\title{\vspace{-.15in}Data-driven Modeling of Rotating Detonation Waves}

\author{ Ariana Mendible$^*$, James Koch$^{**}$, Henning Lange$^\dag$,\\ Steven L. Brunton$^*$, and J. Nathan Kutz$^\dag$\\
{\small \em $^*$Department of Mechanical Engineering, University of Washington, Seattle, WA 98195} \\
{\small \em $^{**}$Oden Institute for Computational \& Engineering Sciences, University of Texas, Austin, TX 78712} \\
{\small \em $^\dag$Department of Applied Mathematics, University of Washington, Seattle, WA 98195} \vspace{-.15in}}
\date{}

\begin{document}

\maketitle

\begin{abstract}
The direct monitoring of a {\em rotating detonation engine} (RDE) combustion chamber has enabled the observation of combustion front dynamics that are composed of a number of co- and/or counter-rotating coherent traveling shock waves whose nonlinear mode-locking behavior exhibit bifurcations and instabilities which are not well understood. Computational fluid dynamics simulations are ubiquitous in characterizing the dynamics of RDE's reactive, compressible flow.  Such simulations are prohibitively expensive when considering multiple engine geometries, different operating conditions, and the long-time dynamics of the mode-locking interactions.  {\em Reduced-order models} (ROMs) provide a critically enabling simulation framework because they exploit low-rank structure in the data to minimize computational cost and allow for rapid parameterized studies and long-time simulations. However, ROMs are inherently limited by translational invariances manifest by the combustion waves present in RDEs.  In this work, we leverage machine learning algorithms to discover moving coordinate frames into which the data is shifted, thus overcoming limitations imposed by the underlying translational invariance of the RDE and allowing for the application of traditional dimensionality reduction techniques. We explore a diverse suite of data-driven ROM strategies for characterizing the complex shock wave dynamics and interactions in the RDE.  Specifically, we employ the dynamic mode decomposition and a deep Koopman embedding to give new modeling insights and understanding of combustion wave interactions in RDEs. 
\end{abstract}

\section{Introduction}
A {\em rotating detonation engine} (RDE) is a novel combustion engine that uses detonative heat release - a nearly constant-volume process - as the dominant mechanism of energy addition to the reactive, compressible fluid flow, contrasting deflagration-based, constant-pressure heat addition typical of aerospace engines. The RDE offers a number of advantages for application in propulsion or land-based power generation, including mechanical simplification, broad operability limits~\cite{Anand2016,Fotia2016}, the potential for increased thermal efficiency~\cite{Nordeen2014,Shao2010}, and the reduction of propellant pumping requirements~\cite{Sousa2017,Rankin2017}.  The operating dynamics of the RDE include co- and counter-rotating coherent combustion wave fronts of varying number which interact to produce a rich set of nonlinear dynamics and instabilities.  Recent modeling efforts have focused on phenomenological models~\citep{Koch2020,Koch2020a} that are capable of reproducing and characterizing the RDE dynamics and bifurcations observed in experiments.  This includes models that characterize the nucleation and formation of combustion pulses, the soliton-like interactions between these combustion fronts, and the fundamental, underlying Hopf bifurcation to periodic modulation of the waves~\citep{Koch2020a}.  The goal of the present work is to characterize the dynamics of the combustion wave front interactions directly from experimental data, specifically with the goal of developing reduced-order models (ROMs) for characterizing the origins of dynamic instabilities in RDEs. We will explore several leading techniques in data-driven optimization (i.e., machine learning) of varying complexity.  

RDE hardware is designed to amplify thermoacoustic instabilities associated with reacting flows in circular and/or periodic geometries. For thrust-producing RDEs, the typical design is an annular combustion chamber, see Figure \ref{fig:intro}a. Fuel and oxidizer are supplied through independent feeds into the head-end of the annulus, where they promptly mix to form a combustible medium.  An ignition source (spark plug) initiates a chemical reaction that quickly and locally releases energy into the fluid. Supposing the geometry of the engine and the rate of heat release allow for a local accumulation of energy (Rayleigh's criterion), sharp gradients in pressure and density (and therefore temperature) form. This creates a feedback loop where chemical kinetics are further accelerated by the increase in temperature, which in turn releases more energy into the fluid. This process saturates once all propellant is locally consumed and combustion halts. However, in the RDE, the sharp gradients in pressure and density form traveling shock waves strong enough to auto-ignite propellant. These shock-reaction structures, or detonation waves, move supersonically about the periodic chamber of the RDE, consuming the newly injected and mixed propellant in its path. The detonations continuously propagate so long as a sufficient amount of mixed propellant exists in its path to overcome dissipative effects (exhaust, for example). A number of experimental RDE programs have detailed the effects of geometry, injection schemes, and fueling conditions (\citep{Dyer2012,Fotia2016,Fotia2017,Walters2019}) on the RDE dynamics.

The detonations follow attractor-like dynamics that are the manifestation underlying multi-scale balance physics of the driven-dissipative RDE~\cite{Koch2020b}.  The RDE is similar in nature to mode-locked lasers~\cite{Haus2000mode,Kutz2006}, where global gain and loss dynamics produce a similar cascading bifurcation diagram of mode-locked states~\cite{li2010geometrical}.  In this context, the mode-locked structures of the RDE are classified as autosolitons, or stably-propagating nonlinear waves where the local physics of nonlinearity, dispersion, gain, and dissipation exactly balance.  These physics are multi-scale in nature: the local fast scale of combustion provides the energy input to generate the mode-locked state, while the slow scales of dissipation and propellant regeneration shape the waveform and dictate the total number of detonation waves. Thus, the global multi-scale balance physics give the detonations their mode-locking properties - not exclusively the frontal dynamics prescribed by classical detonation theory. 

These properties have been experimentally observed at the University of Washington High Enthalpy Flow Laboratory using a gaseous methane-oxygen 76-mm flowpath outside diameter RDE, as described in previous works~\cite{KochPE,Koch2020,Koch2020a}. This experimental apparatus is unique in that the RDE tested is fully modular and that the apparatus exists in a closed system. The modularity of the RDE allows for parametric testing of engine geometries (flowpath lengths and annular gaps) and injectors (varying injection scheme, orifice count, and total injection area) with respect to varied propellant feed rates and stoichiometry. Because the entire apparatus is closed, implied is both the inlet and outlet boundary conditions of the combustor are able to be set. The inlet boundary condition is implicitly set via a desired flow rate and propellant mixture, thereby constraining the manifold pressures. The outlet boundary condition is set via controlling the backpressure of a large (approximately four cubic meters) dump volume. Lastly, the exhaust routing of the engine has allowed for the installation of an optical viewport approximately 2 meters downstream of the exit plane of the combustor. Each experiment consists of four main phases. First, a pre-purge of inert diluent, typically nitrogen, floods the system. Second, the diluent is shut off and propellant begins to flow through the combustor. Third, chemical reactions are triggered, typically via an automotive spark plug or a pre-detonation tube. In a successful experiment, the self-organization of traveling waves occurs and persists so long as propellant is flowing into the combustor. Lastly, the propellant is shut off and diluent is re-introduced into the combustor. For each experiment, a high-speed camera records the duration of the `hot' portion of the run, including the ignition event, the transient mode-locking phase, and steady operation of the combustor. The experiments exhibited in this manuscript are representative of modes of operation and transients observed in this experimental apparatus. These experimental spatiotemporal dynamics are taken from Koch et al.~\cite{KochPE,Koch2020,Koch2020a}.

Computational fluid dynamic (CFD) simulations have been heavily relied upon to diagnose the RDE flowfields. These simulations vary from periodic 2-D `unwrapped' rectangular domains \citep{Schwer2011,Schwer2019,Subramanian2020} to full detailed 3-D engine geometries \citep{Gaillard2017,Sun2017,Sun2018,Lietz2018}. From these simulations, the canonical RDE flowfield is obtained (a cartoon of which is shown in Fig. \ref{fig:intro}a) and relevant metrics can be extracted, such as thrust, specific impulse, available mechanical work, and thermodynamic efficiency \citep{Zhou2012,Nordeen2014}. However, long-time parametric simulations of RDE dynamics is prohibitively expensive since the fastest physics (the detonation front) and the slowest physics (mixing and/or exhaustion) both need to be adequately resolved for proper system behavior.  Thus simulations need to be run for several - if not dozens or hundreds - of cycles (or until the physics of the slowest scales are fully developed). The computational cost of simulations can quickly become prohibitive and it typically requires high-performance computing architectures for even moderate lengths of simulation time.   Consequentially, ROMs have been developed, with varying degrees of success, for recreating the RDE canonical flowfield \citep{Fievisohn2017,Sousa2017a}, predicting thermodynamic trends \citep{Kaemming2017}, predicting application-based propulsive performance \citep{Mizener2017}, or reproducing the dynamics of the waves \citep{DMD_det,Humble2019,Koch2020,Koch2020a}. However, because of the multi-scale nature of the RDE and the intricate interactions of its fundamental physical processes, these modeling efforts are often constrained to geometry, propellant, or mode-specific operating regimes, with the imposition of wave topology or detonation structure. In an alternative approach, experiments  allow us to build ROMs directly from data. 

\begin{figure}
    \centering
    \includegraphics[width = 0.99\textwidth]{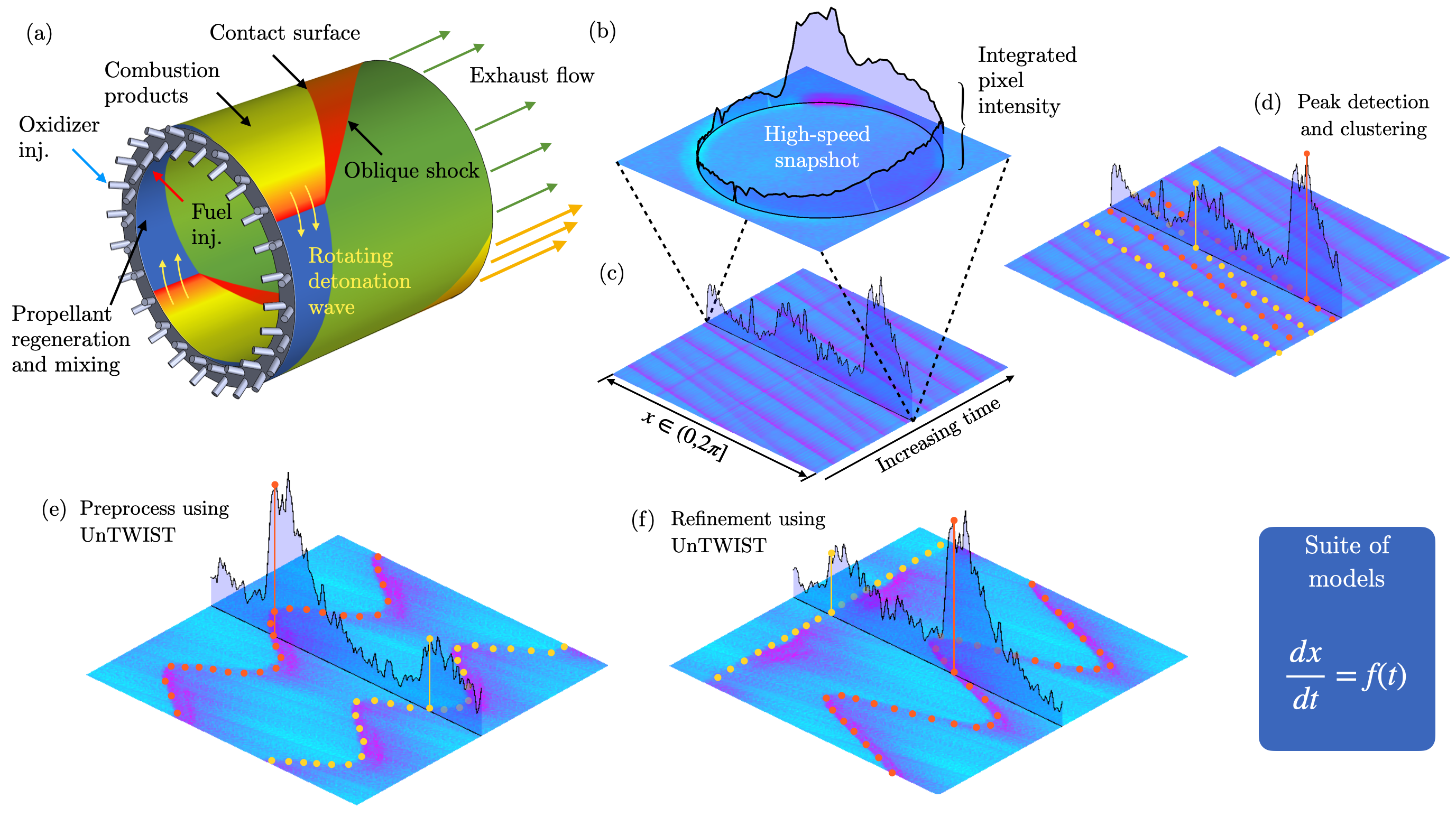}
    \caption{(a) RDE schematic, (b) Schematic of one time slice of video data, viewing down the axis of rotation of the RDE, (c) The same time slice viewed in an $(x,t)$ plot, with each column of the data in time constructed by integrating the pixel intensity along the annulus, (d) demonstration of the peak detection and clustering necessary to model the wave speeds with UnTWIST, (e) a preliminary processing of (c) using the UnTWIST algorithm, (f) a refining processing of (e) with the UnTWIST algorithm, which then becomes the basis for a suite of data-driven models.}
    \label{fig:intro}
\end{figure}

In order to construct ROMs of the combustion-front dynamics, one must first move to a frame of reference of the mode-locked states.  ROMs exploit the intrinsic, low-rank structure of the simulation data in order to create more tractable models for the spatiotemporal evolution dynamics.  Typically ROMs leverage the {\em singular value decomposition} (SVD) to produce a linear dimensionality reduction~\cite{Kutz:2013,Brunton2019book}, whereby a dominant set of correlated modes provide a subspace in which to project the PDE dynamics~\cite{Benner2015siamreview,Taira2017aiaa,taira2019aiaa}.   Typically low-energy modes are truncated, and the governing equations are projected onto the remaining high-energy modes to create an approximate and low-dimensional model. Dimensionality reduction and modal decomposition approaches have been well-studied and are extremely efficient~\cite{antoulas2005approximation,Benner2015siamreview,hesthaven2016certified,quarteroni2015reduced,Taira2017aiaa,taira2019aiaa}.  However, SVD based ROMs are typically compromised by traveling wave physics, which represents an underlying translational invariance.  Thus a growing body of literature is aimed at producing mathematical architectures that are capable of determining the traveling wave frame of reference of the underling wave~\cite{kirby1992reconstructing,Rowley2000physd,Rim2018juq,Reiss2018jsc}.   Mendible et al.~\cite{mendible2020dimensionality} recently developed an unsupervised machine learning procedure for transport-dominated systems characterized by traveling waves.  This method can be applied with or without knowledge of the governing equations, thus providing a robust mathematical architecture for ROMs exhibiting traveling wave phenomenon.   
This algorithmic infrastructure can be used to extract the intrinsic features associated with the RDE front evolution, uncovering a coordinate system where it is possible to obtain low-order models.  
We then leverage a suite of machine learning algorithms to discover the underlying nonlinear dynamics prescribing the ubiquitous RDE front interactions.  Importantly, the methodology is purely data-driven and the ROMs are constructed entirely from detailed experimental observations. 
This work is part of a growing body of literature that is bringing emerging technology in machine learning to bear on problems in fluid mechanics~\cite{duraisamy2019arfm,Brenner2019prf,Brunton2020arfm}.

\section{Detonation Wave Tracking with UnTWIST}

It is widely known that transport phenomenon such as traveling waves impair the effectiveness of traditional dimensionality reduction methods, mainly due to an issue of separation of space and time variables~\cite{Brunton2019book,kirby1992reconstructing,Rowley2000physd,Rim2018juq,Reiss2018jsc}. One approach to resolve this issue is to shift the frame of reference from the laboratory frame to a moving coordinate frame that matches the speed of the traveling waves. Once the traveling quantities have been made stationary in this way, efficient traditional methods such as proper orthogonal decomposition can be utilized for dimensionality reduction.  

In order to build reduced-order models on the RDE data, rife with traveling shock fronts, it is necessary to preprocess it by aligning these traveling waves in time. Here, we employ the unsupervised traveling wave identification with shifting and truncation (UnTWIST)~\cite{mendible2020dimensionality} algorithm to perform this preprocessing step. This method allows for a data-driven and interpretable model for the speeds of the traveling shock fronts, as well as separable low-rank modes, providing an intuitive insight into the physics of the system. A basic overview of using UnTWIST is described here. For further details and a complete algorithm, please refer to~\cite{mendible2020dimensionality}.

Similar to other methods, UnTWIST learns a moving coordinate frame, given by the speed of a traveling wave. This holds the wave of interest stationary, allowing for models to be built for that particular wave. UnTWIST relies on an optimization over a user-input library of potential wave-speed functions to learn this coordinate frame. To execute the optimization on wave profile data $u(x,t)$, two main steps are first performed: (1) ridge detection to learn the location in $(x,t)$ space of the traveling wave fronts or peaks, and (2) spectral clustering to divide the points $(x_i,t_i)$ into groupings for each wave. For example, these two steps are shown in Figure \ref{fig:pks_and_clusters}.

\begin{figure}[t]
\ffigbox
{\begin{subfloatrow}[3] \centering
    \hspace{-0.5cm}
    \setcounter{subfigure}{0}%
    \sidesubfloat[]{\includegraphics[height=0.22\linewidth]{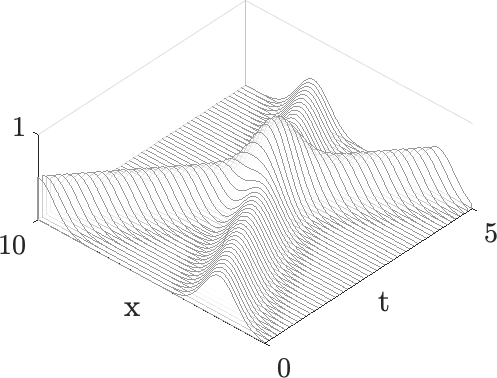}}%
    \hspace{-0.5cm}
    \setcounter{subfigure}{1}%
    \sidesubfloat[]{\includegraphics[height=0.22\linewidth]{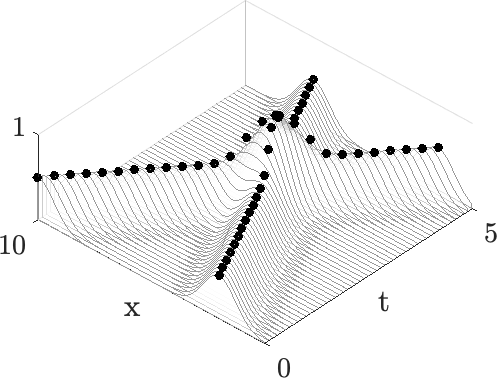}}%
    \hspace{-0.5cm}
    \setcounter{subfigure}{2}%
    \sidesubfloat[]{\includegraphics[height=0.22\linewidth]{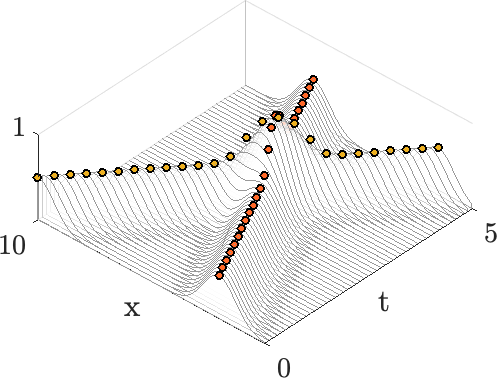}}%
\end{subfloatrow}}
{\addtocounter{figure}{-1}\caption{(a) Example of a traveling wave data set, (b) wave peak points $(x_i,t_i)$ are identified using a ridge detection, shown overlaid with the waves, (c) wave peak points are clustered into wave groups using spectral clustering. Once these points are identified and clustered, a model is fit to each based off of a user-provided library of candidate linear or nonlinear functions.}\label{fig:pks_and_clusters}}
\end{figure}

Once the wave fronts are identified and separated, the data is assembled into the optimization. We construct matrices $\mathbf{X}$ and $\mathbf{T}$ using the $(x,t)$ locations of the wave fronts in $u(x,t)$, where $\mathbf{T}$ contains the values of $t$ evaluated for each function in the user-defined library. The cost function is given by Equation \ref{eq:untwist}

\begin{equation}
	\label{eq:untwist}
	\min_{\mathbf{C}, \mathbf{B}, \mathbf{W} \in \Omega} \frac{1}{2} \mathbf{W \odot \|X-TC}\|^2_2 + \lambda R(\mathbf{B}) + \frac{1}{2\zeta} \|\mathbf{C-B}\|^2_2,
\end{equation}
where $\mathbf{W}$ is the weighting matrix which masks points for clustering into wave groups. $\mathbf{C}$ denotes the coefficients of the models that are discovered for each library term and each wave. We require $\mathbf{C}$ to be sparse by forcing it to be close to an auxiliary  matrix $\mathbf{B}$, which is forced to be sparse via a regularizing function $\mathbf{R}$. The hyperparameter $\lambda$ is chosen to enforce a sparse solution, i.e. only a small number of library functions are used to model each wave. The hyperparameter $\zeta$ is chosen to enforce relaxation of the sparsity of $\mathbf{C}$. This optimization presents a large search space over multiple parameters, and is not guaranteed to be convex. Sparse relaxed regularized regression (SR3)~\cite{zheng2018unified} is used to minimize the cost function because of its ability to handle non-convexity and its computational efficiency compared to similar sparsity-promoting optimization schemes.

Once the model coefficients $\mathbf{C}$ are learned, they can be used to shift each time segment of data in order to align the data into one wave group's moving coordinate frame. The mask matrix $\mathbf{W}$ allows for easy separation of the wave fronts for this alignment. Once the data is aligned into the new coordinate frame, traditional dimensionality reduction methods can easily be applied and can be expected to reveal extremely low-rank modes for the `straightened' wave or wave group. 

\subsection{UnTWIST Applied to RDE Data} 
An example of the UnTWIST algorithm applied to snapshots of RDE data can be seen in Figure \ref{fig:methods_example}.
UnTWIST was applied in two steps for each data set presented. The models we build are based on time series that are long relative to the spatial dimension-- with between 1,000 and 10,000 time steps relative to 180 or 360 spatial points. The snapshots also contain wave fronts that travel on a fast time scale relative to the slow time scale of the relevant dynamics, see Figure \ref{fig:methods_example}a, necessitating an extreme shift in order to shift into a straightened wave coordinate frame. Because of the fast-moving fronts and long time series, the UnTWIST algorithm was applied in two steps-- a preprocessing step, and a refining step, with different inputs for each. 

For the preprocessing step, only 10 time steps of the data are considered, as seen in Figure \ref{fig:methods_example}b. Using the UnTWIST algorithm and the identified wave peak points as shown in Figure \ref{fig:methods_example}c, we obtain a simple linear model to be fit to each wave, Figure \ref{fig:methods_example}d. A single linear-speed shift, the average of the speed models, is applied to original data for the entire time series, and gives Figure \ref{fig:methods_example}e. This first shift reveals critical underlying dynamics of the shock wave interactions. 

\begin{figure}[t]
\ffigbox
{\begin{subfloatrow}[1]
\sidesubfloat[]{\includegraphics[height=0.22\linewidth]{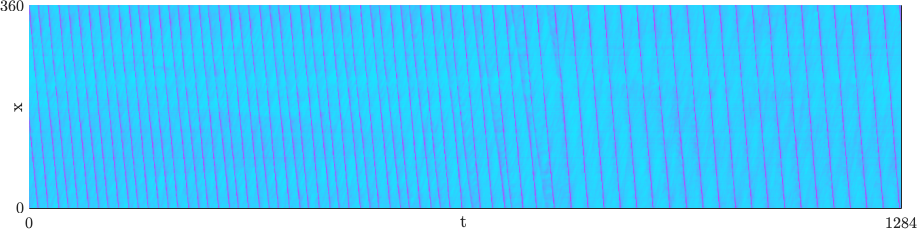}\label{trans:orig}}%
\end{subfloatrow}\par\bigskip
\begin{subfloatrow}[3] \centering
\setcounter{subfigure}{1}%
\sidesubfloat[]{\includegraphics[height=0.22\linewidth]{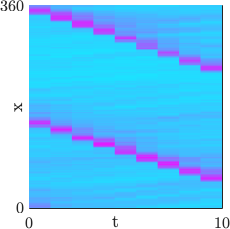}\label{trans:10sec}}%
\setcounter{subfigure}{2}%
\sidesubfloat[]{\includegraphics[height=0.22\linewidth]{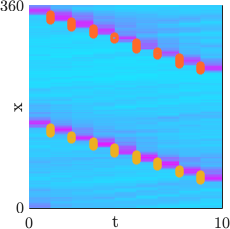}\label{trans:pts}}%
\setcounter{subfigure}{3}%
\sidesubfloat[]{\includegraphics[height=0.22\linewidth]{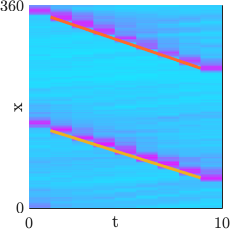}\label{trans:models}}%
\end{subfloatrow}\par\bigskip
\begin{subfloatrow}[1]
\sidesubfloat[]{\includegraphics[height=0.22\linewidth]{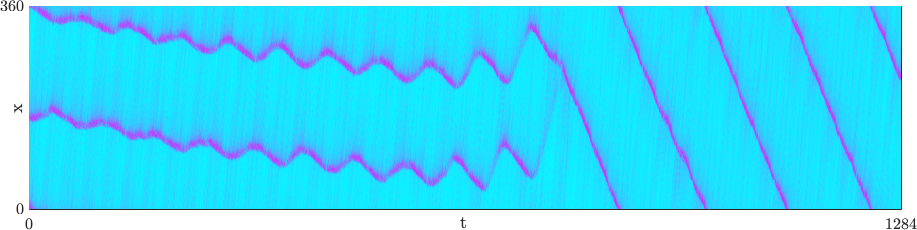}\label{trans:shift}}%
\end{subfloatrow}}
{\addtocounter{figure}{-3}\caption{Preprocessing step using UnTWIST: (a) Example of an original data set, presenting fast-moving shock fronts and a long time series relative to the spatial dimension, (b) A 10 time step segment of (a) showing approximately linear-speed shock front propagation, (c) wave peak points identified and separated, (d) linear models of the shock wave speeds, (e) Data from (a) shifted into the average wave speed, dictated by the models identified in (d), which reveals the relevant interactions between the shock waves.}\label{fig:methods_example}}
\end{figure}

After the first shift is performed, a second refining shift can be used with a more diverse library of potential wave speed models to completely align the data for building low-rank models. The refining shift is performed similarly to the first shift, but we now include potential wave speed functions such as sinusoids, polynomials, exponentials, and nonlinear combinations of these terms. For the example shown in Figure \ref{fig:methods_example2}, sinusoids and exponentials were included in the candidate function library. Once the models were learned, the peak points were refined and the data was shifted into the a refined coordinate frame, one for each shock wave.

The outcome of the second shift can be seen in Figure \ref{fig:methods_example2}. Each coordinate frame allows one shock wave to appear stationary at a time. the shifted data is now aligned in a manner that is amenable to traditional dimensionality reduction methods, such as POD. An example of the first mode of a robust dimensionality reduction~\cite{Candes:2011} of the shifted data is shown in Figures \ref{fig:methods_example2}b and \ref{fig:methods_example2}d, compared to the first mode of the POD of the laboratory frame original data. Here, we only explore the first segment of the time series, before the bifurcation point. 

\begin{figure}[tb]
\ffigbox
{
\begin{subfloatrow}[2]
\setcounter{subfigure}{0}%
\sidesubfloat[]{\includegraphics[width=0.53\linewidth]{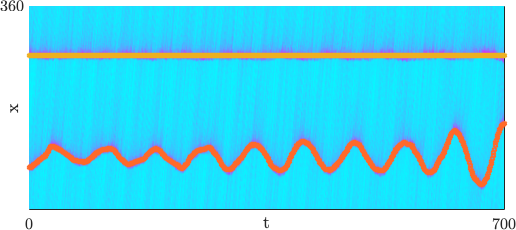}\label{trans:shift1}}%
\setcounter{subfigure}{1}%
\sidesubfloat[]{\includegraphics[width=0.35\linewidth]{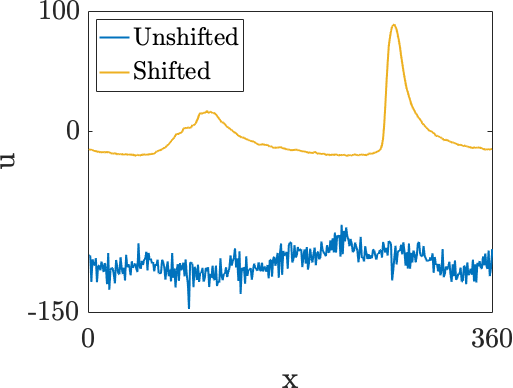}\label{trans:mode1}}%
\end{subfloatrow}\par\bigskip
\begin{subfloatrow}[2]
\setcounter{subfigure}{2}%
\sidesubfloat[]{\includegraphics[width=0.53\linewidth]{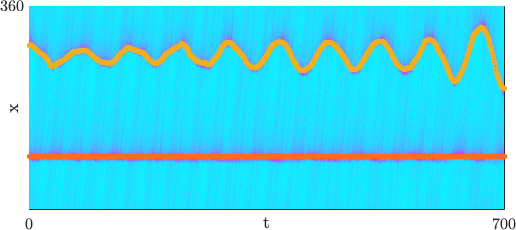}\label{trans:shift2}}%
\setcounter{subfigure}{3}%
\sidesubfloat[]{\includegraphics[width=0.35\linewidth]{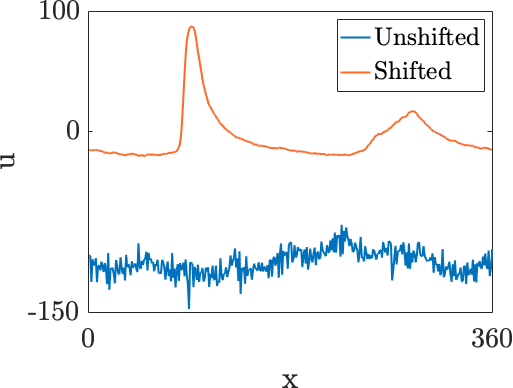}\label{trans:mode2}}%
\end{subfloatrow}
}
{\addtocounter{figure}{-2}\caption{Refinement step using UnTWIST: (a) Data shifted into coordinate frame of the first (yellow) wave. (b) Corresponding first mode of a robust dimensionality reduction of the shifted data in (a) compared to the first mode of the original data. The shifted mode shows a clear shock front where the first wave is straightened, and a smaller artifact where the second wave exists, whereas the original data mode reveals no interpretable structures within the data. (c) Data shifted into the coordinate frame of the second (orange) wave. (d) First mode from dimensionality reduction of the second wave frame, showing the same shock front shape in the correct position, and artifact of the first wave.}\label{fig:methods_example2}}
\end{figure}

This provides an example of how UnTWIST is used on a particular data set in order to align the traveling waves to uncover low-rank representations of their wave fronts. The same steps have been used to process various data sets. While UnTWIST can align these wave fronts into more amenable coordinate frames for dimensionality reduction of the wave field as a whole, it is also of great interest to study the linear and nonlinear interactions between wave fronts in RDEs. Using similarly aligned data and the wave speeds and locations throughout the time series, we explore models of the shock wave interactions in the following section. 

\section{Data-Driven Models of Rotating Detonation Front Dynamics}

The ability to automate the discovery of a moving coordinate system pinned to a shock front allows for a wide range of reduced-order modeling possibilities.  In what follows, we develop a number of data-driven modeling strategies that reduce the dynamics to simple models that characterize the observed interactions of the rotating detonation waves.  Experiments show that these interactions can range from simple linear dynamics to more complicated nonlinear dynamical interactions.  Our diversity of methods allow us to characterize the full range of observed data.

\subsection{Linear Dynamics: Dynamic Mode Decomposition} 

The dynamic mode decomposition (DMD)~\cite{Schmid2010jfm,Rowley2009jfm,Tu2014jcd,Kutz2016book} is an alternative to the the {\em proper orthogonal decomposition} (POD) reduction typically used in ROMs.  It not only correlates spatial activity, but also enforces that various low-rank spatial modes be correlated in time, essentially merging the favorable aspects of POD in space and the Fourier transform in time.  Thus in addition to performing a low-rank SVD approximation, it further performs an eigendecomposition on a best-fit linear operator that advances measurements forward in time in the computed subspaces in order to extract critical temporal features.  
The DMD algorithm is a least-square regression. In its simplest form~\cite{Tu2014jcd}, one can consider two sets of measurement data
\begin{equation}
  {\bf X} = \begin{bmatrix}
\vline & \vline & & \vline \\
{\bf u}_1 & {\bf u}_2 & \cdots & {\bf u}_{m-1}\\
\vline & \vline & & \vline
\end{bmatrix} \,\,\,\,\, \mbox{and} \,\,\,\,\, 
{\bf X}' = \begin{bmatrix}
\vline & \vline & & \vline \\
{\bf u}'_1 & {\bf u}'_2 & \cdots & {\bf u}'_{m-1}\\
\vline & \vline & & \vline
\end{bmatrix}
\end{equation}
where the primed data is advanced $\Delta t$ into the future compared to its unprimed counterpart.   Exact DMD computes the leading eigendecomposition of the best-fit linear operator ${\bf A}$ relating the data
\begin{equation}
  {\bf A} = {\bf X}' {\bf X}^\dag.
  \label{eq:DMD}
\end{equation}
where $\dag$ represents the Moore-Penrose pseudo-inverse.  This gives a least-square fit to the best linear model fitting the data whose solution is
\begin{equation}
   {\bf u}_k=  \sum_{j=1}^n  \boldsymbol{\phi}_j \lambda_j^k b_j=\boldsymbol{\Phi} \boldsymbol{\Lambda}^k\mathbf{b}\,
   \label{eq:omegaj}
\end{equation}
where $\boldsymbol{\phi}_j$ and $\lambda_j$ are the eigenvectors  and eigenvalues of the matrix ${\bf A}$, and the coefficients $b_j$ are the coordinates of the initial condition ${\bf u}_0$ in the eigenvector basis. The eigenvalues $\lambda$ of ${\bf A}$ determine the temporal dynamics of the system.  It is often convenient to convert these eigenvalues to continuous time, $\omega = \log(\lambda)/\Delta t$, so the real parts of the eigenvalues $\omega$ determine growth and decay of the solution, and the imaginary parts determine oscillatory behaviors and their corresponding frequencies.  The eigenvalues and eigenvectors are critically enabling for producing interpretable diagnostic features of the dynamics.  

For the present purposes, the DMD algorithm can be simplified even further by simply fitting the linear dynamics models to the following linear oscillator with damping.  The oscillator takes the continuous-time DMD form
\begin{equation}
    \overbrace{
    \begin{bmatrix}
    \dot{\mathbf{x_1}} \\
    \dot{\mathbf{x_2}}
    \end{bmatrix}}^
    {\dot{{\bf x}}}
     = 
     \overbrace{
     \begin{bmatrix}
     0 & 1 \\
     \alpha & 0 
     \end{bmatrix}}^{\mathbf{A}}
     \overbrace{
     \begin{bmatrix}
     \mathbf{x_1}\\ \mathbf{x_2}
     \end{bmatrix}}^{\mathbf x}
\end{equation}
where $\mathbf{x_1}$ and $\mathbf{x_2}$ represent the normalized distance and velocity between the waves, respectively and $\alpha$ is frequency of oscillation found by least-square regression.  Damping can be added to capture the transient growth observed.  
We solve directly for $\mathbf{A}$ and then find its eigenvalues $\mathbf{D}$ and eigenvectors $\mathbf{V}$. Using the initial conditions $\mathbf{x_1}(0)$ and $\mathbf{x_2}(0)$ directly from the data, we obtain the simple linear model of the distance between the wave fronts as 
\begin{equation}\label{eq:msd}
    x_1 (t) = 18.38 \exp{(0.0011 t)} \cos{(0.0837t)}. 
\end{equation}
Figure \ref{fig:dmd-model} shows the true distance between the shock fronts compared to the linear model in \ref{eq:msd}. 
This simplified DMD model circumvents many of the biases introduced by standard DMD algorithms~\cite{Schmid2010jfm,Rowley2009jfm,Tu2014jcd,Kutz2016book} by directly fitting an exponential solution (optimized DMD)~\cite{askham2018variable}. 

\begin{figure}[t]
    \centering
    \includegraphics[width = 0.5\linewidth]{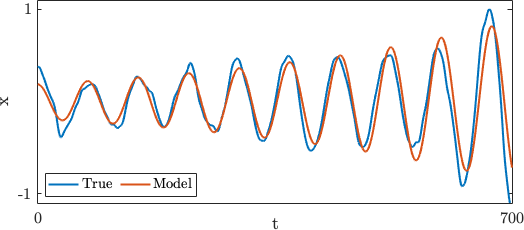}
    \caption{True distance between the shock wave fronts compared to the mass-spring-damper oscillator model given in Eq.~\ref{eq:msd}.}
    \label{fig:dmd-model}
\end{figure}

\subsection{Nonlinear Dynamics: Lotka-Volterra Model} 
Although there are many interacting wave dynamics that appear to be described by simple harmonic motion, i.e. linear oscillators well-captured by DMD, the RDE also produces dynamics that are strongly nonlinear in nature.
Figure \ref{fig:291_data}a presents three RDE shock fronts interacting in an oscillatory manner. When shifted into the coordinate frame of the top wave, shown in \ref{fig:291_data}b, it is clear that the middle and lower wave exhibit sharp, periodic changes in wave speed with respect to the top wave.  Such oscillations are beyond a simple linear description. 

\begin{figure}[b]
\ffigbox
{\begin{subfloatrow}[1]
\setcounter{subfigure}{0}%
\sidesubfloat[]{\includegraphics[width=0.75\linewidth]{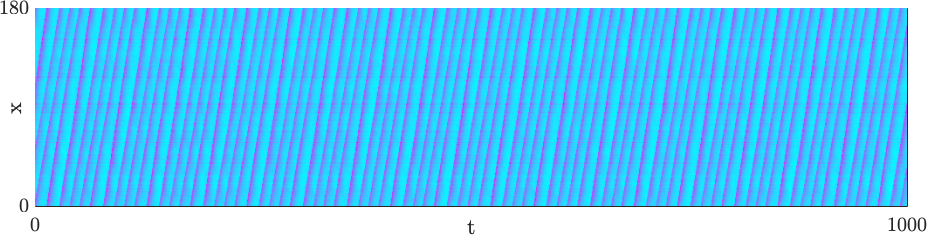}}%
\end{subfloatrow}\par\bigskip
\begin{subfloatrow}[1]
\sidesubfloat[]{\includegraphics[width=0.75\linewidth]{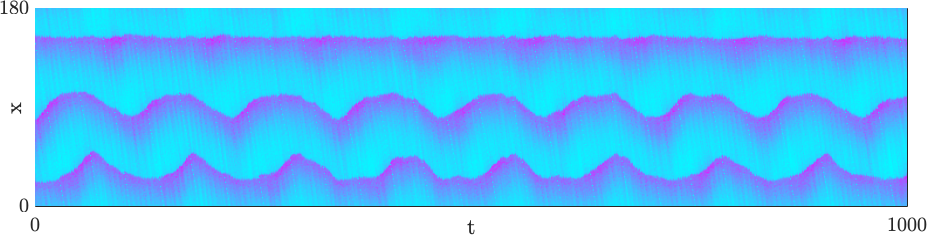}}
\end{subfloatrow}\par\bigskip}
{\addtocounter{figure}{-2}\caption{Figure (a) an example data set in the laboratory frame. (b) Same data set shifted into the coordinate frame learned by UnTWIST using linear models for wave speed for both the preprocessing and refinement steps. }\label{fig:291_data}}
\end{figure}

\begin{figure}
    \centering
    \includegraphics[width=0.7\linewidth]{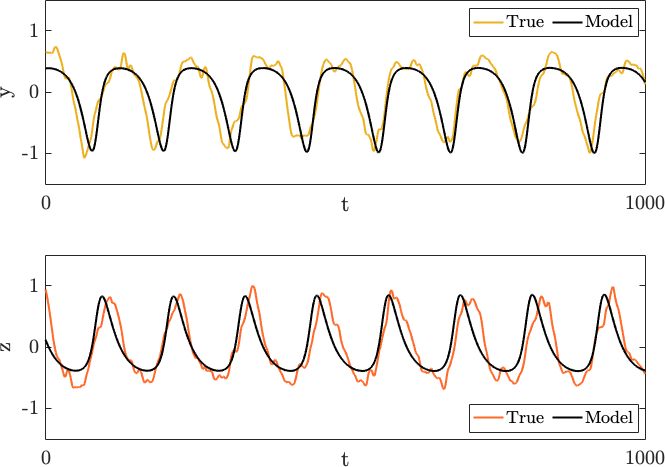}
    \caption{Lotka-Volterra model for the two traveling waves in Figure \ref{fig:291_data}. Yellow and orange indicate the true peak locations $y$ and $z$ of the two wave fronts of the data set in the shifted frame, and black shows the Lotka-Volterra model. Parameters of $\alpha = 0.07$, $\beta = 0.13$, $\delta = 0.10$ and $\gamma = 0.05$ yield a well-fitting model for the two waves which compete for resources.}
    \label{fig:lotka-volterra}
\end{figure}

We explore this nonlinear dynamics by fitting a Lotka-Volterra model to the peak locations of the middle and lower waves. These locations are denoted as $y$ and $z$. The Lotka-Volterra equations, also known as the predator-prey equations, are a coupled pair of nonlinear equations often used to model population changes in two species and are given by:
\begin{align}\label{eq:lotka-volterra}
    \frac{dy}{dt} &= \alpha y - \beta yz \\
    \frac{dz}{dt} & = \delta yz - \gamma z, 
\end{align}
where $\alpha$, $\beta$, $\delta$, and $\gamma$ are positive real parameters controlling the growth and decay of $y$ (prey), and growth and decay of $z$ (predators), respectively. We take $y$ and $z$ to be the locations of the two traveling waves from Figure \ref{fig:291_data}b, using the negative of the middle wave to orient the sharp peaks in the positive direction. Parameters for the best-fit model were chosen via a parameter sweep, and were found to be $\alpha = 0.07$, $\beta = 0.13$, $\delta = 0.10$ and $\gamma = 0.05$. The first 500 time steps were used as a training set. The error was computed over these time steps as the Frobenius norm of the difference between model and true $[y,z]$. The resulting model, with extrapolation out to 1000 time steps can be seen in Figure \ref{fig:lotka-volterra}.

This model proves to be a good fit for the periodic nonlinear dynamics of wave interactions, with the frequency matching through to the end of the test data set. Additionally, the Lotka-Volterra model intuitively fits the nature of the data. In the RDE, each detonation wave competes for the shared fuel resource as it travels through the combustion chamber. One wave's acceleration consumes more fuel along the annulus, and therefore leads to the other waves' deceleration, similar to how a decrease in population of prey leads to a decrease in predator population. The Lotka-Volterra model is only one of a potential wealth of interpretable nonlinear models that may describe and give insight into the physics governing the RDE. 
An interesting avenue of future work is to automate the identification of nonlinear dynamics, for example using the sparse identification of nonlinear dynamics (SINDy)~\cite{Brunton2016pnas} algorithm.  SINDy has been widely applied to identify reduced-order models for fluid systems~\cite{Loiseau2017jfm,Loiseau2018jfm}, and is a promising candidate for obtaining low-order models of RDE dynamics.  

\subsection{Deep Koopman}

In the following, ROMs based on Koopman theory are explored that do not rely on an explicit dimensionality reduction technique and can therefore avoid translational invariances. Koopman theory postulates that any nonlinear dynamical system can be lifted by the means of a nonlinear and time-invariant functional, oftentimes referred to as observables, into a space where its time evolution can be described by linear methods. It was first introduced in the seminal 1931 paper for Hamiltonian systems~\citep{koopman1931hamiltonian} but later generalized to continuous-spectrum systems \citep{koopman1932dynamical}. Even then it was of considerable importance as a building block for advances in ergodic theory~\citep{birkhoff1931proof,neumann1932proof,birkhoff1932recent,neumann1932physical,moore2015ergodic}. Koopman theory has experience renewed interest in the past two decades~\cite{mezic2005spectral,budivsic2012applied,Mezic2013arfm}.

Let $u(t)$ be the collected measurements at time $t$ and $\psi$ be the time-invariant observable functional, Koopman theory dictates that there a linear operator $\mathcal{K}$ always exists such that:
\begin{align*}
    \mathcal{K}\psi(u(t)) = \psi(u(t+1)).
\end{align*}
$\mathcal{K}$ is usually referred to as the Koopman operator and might in practice be infinite dimensional as for example in chaotic systems.

Recent applied research has focused on algorithmic approaches to estimate the Koopman operator from measurement data. Early approaches relied on auto-encoder structures~\citep{otto2019linearly,lusch2018deep,yeung2019learning,wehmeyer2018time,takeishi2017learning} and solved an optimization objective, usually by means of gradient descent, usually composed of terms that encourage linearity in `Koopman space' and reconstruction performance. These approaches were extended in various ways. For example, Bayesian Neural Networks as encoders were utilized to extend Koopman theory to the probabilistic setting~\citep{pan2019physics}. Furthermore Champion et al.~\cite{champion2019data} relaxed the linearity requirement and allowed for sparse dynamics in the latent space. Because linearity in `Koopman space' is part of the optimization objective, these approaches are usually only approximately and locally linear, which in turn impedes their ability to predict far into the future. A more recent approach called Koopman forecast~\cite{lange2020fourier} is linear in `Koopman space' by construction and does not require training an encoder network. In order to overcome its nonlinear and non-convex objective, the Koopman Forecast algorithm employs gradient descent in conjunction with the Fast Fourier Transform. In the following, variants of the Koopman Forecast algorithm are introduced and their efficacy in modeling rotating detonation waves is evaluated.

\subsubsection{Koopman Forecast}

The Koopman Forecast algorithm provides the tools to approximate the Koopman operator from data that is assumed to be quasi-periodic. The quasi-periodicity assumption in turn restricts the Koopman operator to have purely imaginary eigenvalues, i.e. the Koopman operator describes a stable linear dynamical system.

Note that for any linear system $y(t)$ with purely imaginary eigenvalues the following holds:

\begin{align*}
    y(t) \propto \begin{bmatrix}
                           \cos(\omega_1 t)\\
                           \vdots \\
                           \cos(\omega_{n} t)\\
                           \sin(\omega_1 t)\\
                           \vdots \\
                           \sin(\omega_{n} t)\\
                    \end{bmatrix} := \Omega(\boldsymbol{\omega} t).
\end{align*}

Because of this, the Koopman Forecast algorithm solves the following optimization problem:
\begin{align*}
    E(\vec{\omega}, \Theta) = \sum_{t=1}^T ||u(t) - f_\Theta(\Omega(\boldsymbol{\omega} t)) ||_2^2.
\end{align*}

In this case, $f_\Theta$ is some nonlinear function parameterized by $\Theta$, for example a Neural Network. Thus, colloquially speaking, the Koopman Forecast algorithm fits a Neural Network driven by a linear oscillator to data. Because of the nonlinearity of $f$, the objective $E$ is notoriously difficult to solve for $\boldsymbol{\omega}$. However, as laid out in~\cite{lange2020fourier}, by exploiting periodicities in temporally local loss functions in conjunction with coordinate descent, globally optimal values in the direction of $\omega_i$ can be obtained. Specifically, for every $t$, the temporally local loss function $||u(t) - f_\Theta(\Omega(\boldsymbol{\omega} t)) ||_2^2$ is periodic in $\frac{2\pi}{t}$. Thus, it is sufficient to sample each temporally local loss function within its first period in order to reconstruct $E$. The reader is referred to~\cite{lange2020fourier} for a more thorough discussion of this technique.

\subsubsection{Temporal Koopman with spatial decoder}

\begin{figure}
    \centering
    \includegraphics[width=0.68\linewidth]{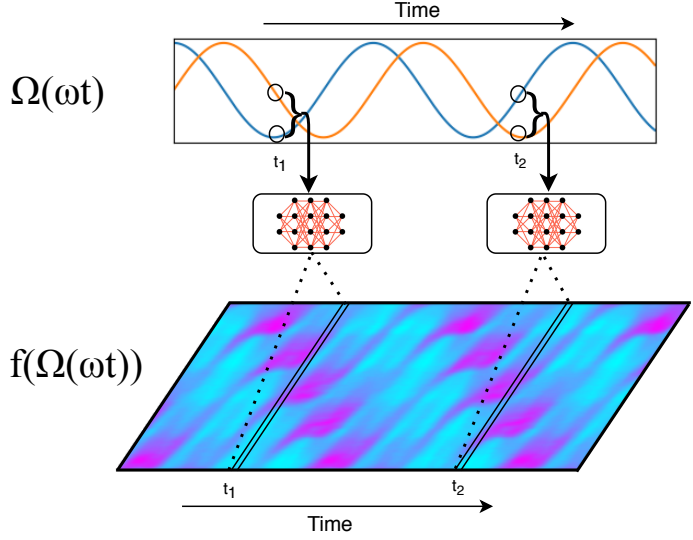}
    \caption{A graphical depiction of the Koopman Forecast algorithm with one underlying frequency. Colloquially speaking, the Koopman Forecast algorithm resembles a Neural Network driven by a linear oscillator.}
    \label{fig:kf_diag}
\end{figure}

In the following, we will show how the Koopman Forecast with a spatial decoder can alleviate the problem of translational invariance. For this, a spatial decoder function $f_\Theta$ is devised that converts time into space. As $f_\Theta$, we choose a fully connected feed-forward Neural Network with the following topology: $2 \rightarrow 32 \rightarrow 32 \rightarrow 180$ and $\tanh$ nonlinearity in intermediate layers. The output dimensionality is 180 because space is sampled at 180 locations. The input dimensionality is 2 because we assume the system to be driven by a single frequency. Figure \ref{fig:kf_diag} shows graphically the setup for the experiment. Note that $f_\Theta$ can learn and preserve the nonlinear interactions between waves but because it is a Neural Network, it is difficult to extract interpretable information about the nature of the nonlinear interactions. Figure \ref{fig:008_direct} shows the spatiotemporal prediction of the system into the future (b) against the true data (a), with Figure \ref{fig:008_direct}c comparing a single time slice. The algorithm correctly explains approximately 75\% of the variance. Considering the measurements are obtained by experimentation and therefore exhibit a considerable amount of noise, the Koopman Forecast algorithm seems to perform well.

\begin{figure}[tb]
\ffigbox
{\begin{subfloatrow}[1]
\setcounter{subfigure}{0}%
\sidesubfloat[]{\includegraphics[width=0.75\linewidth]{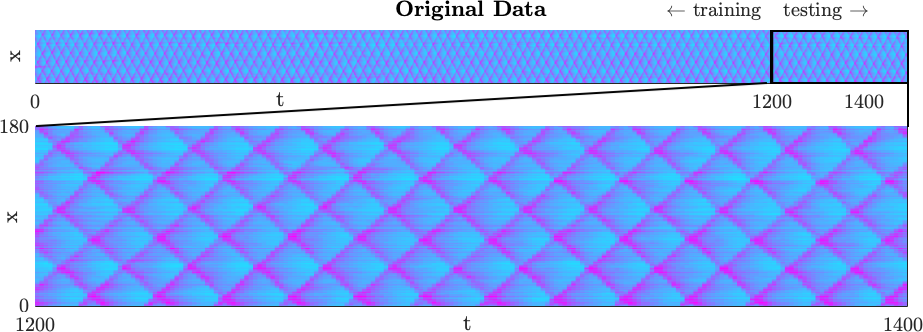}}%
\end{subfloatrow}\par\bigskip
\begin{subfloatrow}[1]
\sidesubfloat[]{\includegraphics[width=0.75\linewidth]{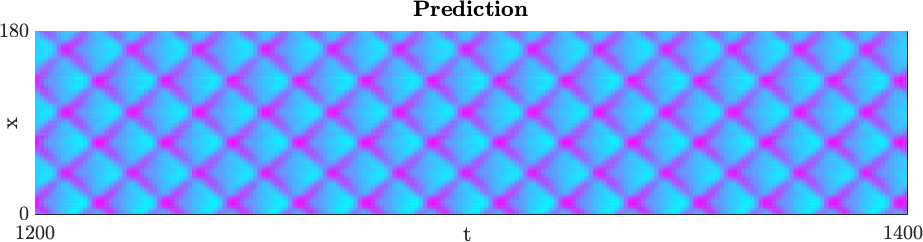}}
\end{subfloatrow}\par\bigskip
\begin{subfloatrow}[1]
\sidesubfloat[]{\includegraphics[width=0.45\linewidth]{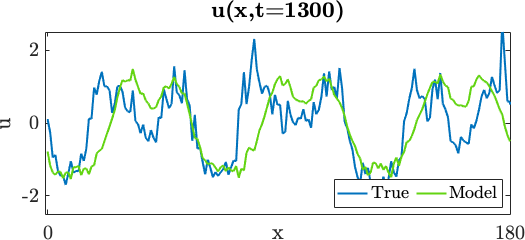}}
\end{subfloatrow}\par\bigskip}
{\addtocounter{figure}{-3}\caption{(a) Original data in the laboratory coordinate frame, showing training data from times $t=[0, 1200)$ and testing data from times $t=[1200, 1400]$, and a broader view of the testing data set. (b) Prediction of the Koopman Forecast algorithm over the same testing time. (c) Comparison of the prediction of the wave shape using the Koopman Forecast algorithm to the ground truth wave shape at an example time step of $t=1300$.}\label{fig:008_direct}}
\end{figure}

\begin{figure}[t]
\ffigbox
{\begin{subfloatrow}[1]
\setcounter{subfigure}{0}%
\sidesubfloat[]{\includegraphics[width=0.75\linewidth]{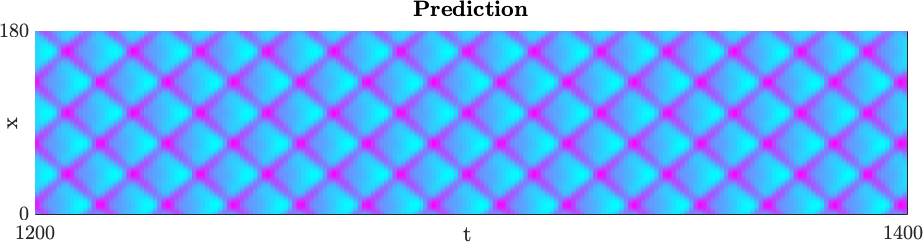}}
\end{subfloatrow}\par\bigskip
\begin{subfloatrow}[1]
\sidesubfloat[]{\includegraphics[width=0.75\linewidth]{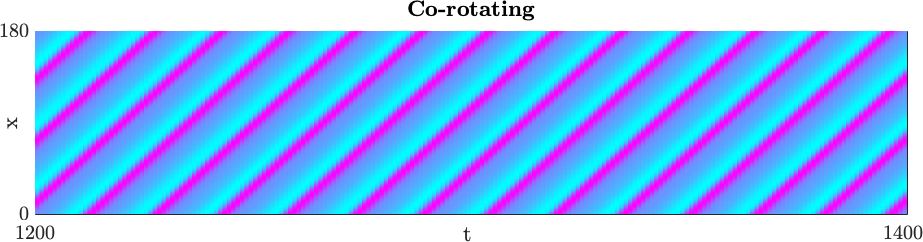}}%
\end{subfloatrow}\par\bigskip
\begin{subfloatrow}[1]
\sidesubfloat[]{\includegraphics[width=0.75\linewidth]{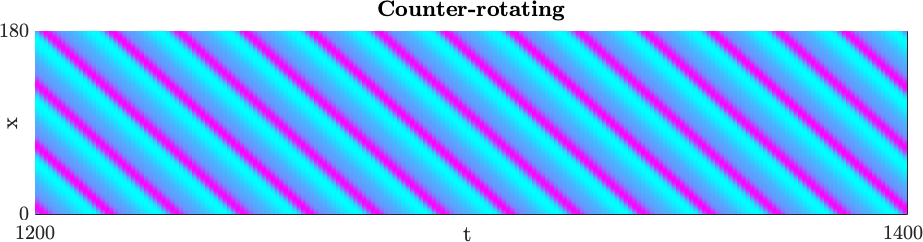}}
\end{subfloatrow}\par\bigskip
\begin{subfloatrow}[1]
\sidesubfloat[]{\includegraphics[width=0.45\linewidth]{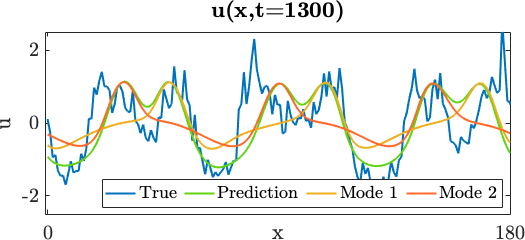}}
\end{subfloatrow}\par\bigskip}
{\addtocounter{figure}{-4}\caption{(a) Shows the prediction of the modal Koopman model on the same data from Figure \ref{fig:008_direct}a, (b) and (c) show the inferred modes respectively. Note that colloquially, (a) = (b) + (c). (d) shows the inferred modes at the same time step and the aggregate prediction (Mode 1 + Mode 2) in comparison to the ground truth measurement data at an example time $t=1300$.}\label{fig:008_modes}}
\end{figure}

\begin{figure}[t]
\ffigbox
{\begin{subfloatrow}[1]
\setcounter{subfigure}{0}%
\sidesubfloat[]{\includegraphics[height=0.19\linewidth]{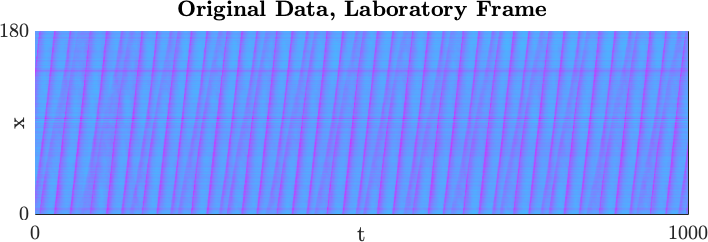}}%
\sidesubfloat[]{\includegraphics[height=0.19\linewidth]{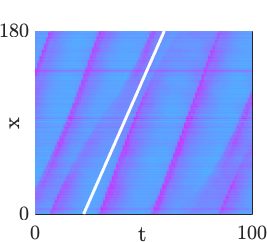}}
\end{subfloatrow}\par\bigskip
\begin{subfloatrow}[1]
\sidesubfloat[]{\includegraphics[height=0.19\linewidth]{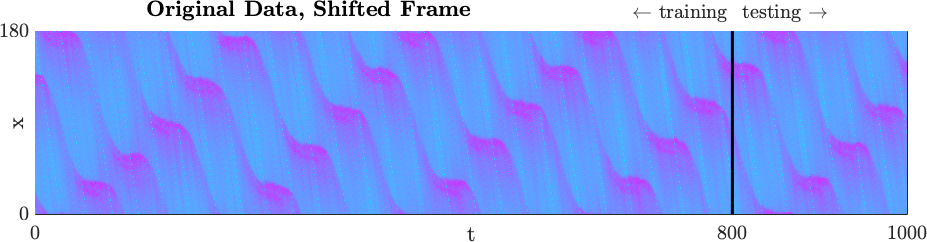}}
\end{subfloatrow}\par\bigskip
\begin{subfloatrow}[1]
\sidesubfloat[]{\includegraphics[height=0.19\linewidth]{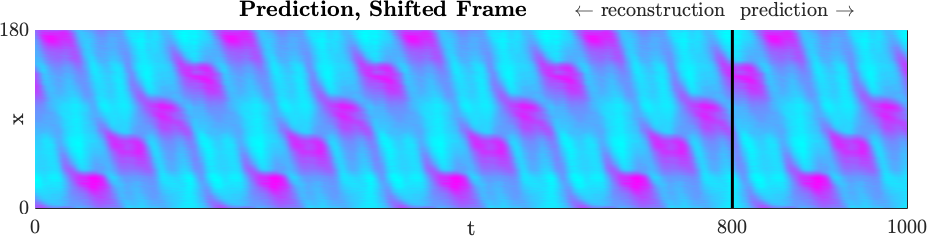}}%
\end{subfloatrow}\par\bigskip
\begin{subfloatrow}[1]
\sidesubfloat[]{\includegraphics[height=0.19\linewidth]{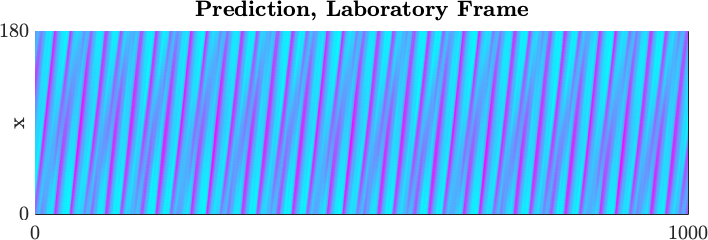}}%
\sidesubfloat[]{\includegraphics[height=0.19\linewidth]{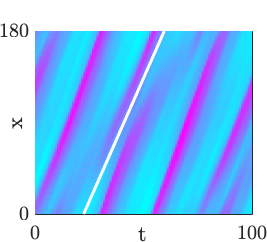}}
\end{subfloatrow}\par\bigskip
\begin{subfloatrow}[1]
\sidesubfloat[]{\includegraphics[height=0.19\linewidth]{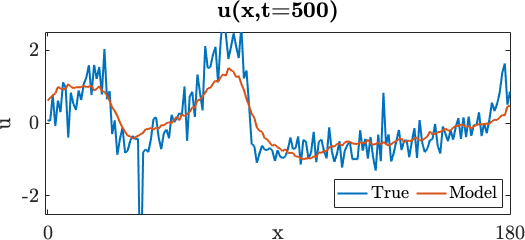}}%
\end{subfloatrow}\par\bigskip}
{\addtocounter{figure}{-5}\vspace{-.2in}\caption{A comparison of the true measurements to the prediction of the Koopman Forecast algorithm. (a) the full time series of original data in the laboratory frame. (b) highlights a short time series showcasing the nonlinear interactions between wave fronts, as they do not run parallel to the white guideline. (c) data shifted into the coordinate frame learned through UnTWIST. (d) presents the prediction of the same data using the Koopman Forecast algorithm; (e) and (f) show the Koopman Forecast shifted back into the laboratory frame, comparing to (a) and (b) respectively. It is clear in (f) that the nonlinear interaction between the waves is well-preserved. (g) compares the learned model to the true data at one time slice, $t=500$, showcasing the robustness to noise and the de-noising effect of the model.}\label{fig:387}} 
\end{figure}

\subsubsection{Spatiotemporal Koopman}

The Koopman Forecast model paired with a spatial decoder seems to explain the data reasonably well. However, because $f_\Theta$ is a Neural Network, it is hard for practitioners to understand what the algorithm has learned, i.e. the model gives little insight apart from the extracted frequency. In the following a more idealized model will be introduced that is less flexible but also more interpretable. The assumptions of the model are the following: We assume that $N$ modes interact linearly and that both modes travel at a constant speed. Specifically, we model the data as a superposition of $N$ rotating or shifted modes $m_i$.

Let $u(x,t)$ denote the wave height at position $x$ and time $t$. Mathematically speaking, we assume the following:
\begin{align*}
    u(x,t) = \sum_i^N m_i(f_i(t)+x),
\end{align*}
where $f_i(t)$ is a time-dependent function that models the offset of mode $i$ at time $t$. In a next step, we incorporate the knowledge that spatial boundary conditions are periodic. Assume that $u \in [0, K]$. We can find a periodic parameterization of $m_i$ in the following way:
\begin{align*}
    m_i(x) = g_{\theta_i}\left(\begin{bmatrix}
           \sin(\frac{2 \pi}{K} x) \\
           \cos(\frac{2 \pi}{K} x) 
         \end{bmatrix}\right)
\end{align*}
with $\theta_i$ being model parameters of the $i$th mode (e.g. weights of a Neural Network). Note that this parameterization is periodic in $K$ because:
\begin{align*}
    m_i(x + K) = g_{\theta_i}\left(\begin{bmatrix}
           \sin(\frac{2 \pi}{K} x + 2\pi) \\
           \cos(\frac{2 \pi}{K} x + 2\pi) 
         \end{bmatrix}\right) = m_i(x).
\end{align*}

If we now assume that the offset of the $i$th mode increases at a constant speed $\omega_i$, we can rewrite $u(x,t)$ in the following way:
\begin{align*}
    u(x,t) = \sum_i^N  g_{\theta_i}\left(\begin{bmatrix}
           \sin(\frac{2 \pi}{K} x + \omega_i t) \\
           \cos(\frac{2 \pi}{K} x + \omega_i t) 
         \end{bmatrix}\right).
\end{align*}
Let $\Theta = \{ \theta_i \}_{i=1}^N$. Fitting $u(x,t)$ to measured data $y(x,t)$ requires solving:
\begin{align*}
    E(\vec{\omega}, \Theta) = \sum_t \sum_x (y(x,t) - u(x,t)))^2.
\end{align*}
Solving this optimization objective for $\omega$ is, again, notoriously difficult as it is not only non-convex but also non-linear. However, note that for every $x$ and $t$, $(y(x,t) - u(x,t)))^2$ is, again, periodic in ${2\pi}/{t}$. Therefore, an analogous strategy to the Koopman Forecast algorithm can be employed to solve for $\boldsymbol{\omega}$. 

Figure \ref{fig:008_modes} shows the results when $N=2$, i.e. when two additive modes are assumed. The modal Koopman model is considerably stiffer as it only explains 65\% of the variance in comparison to the Koopman Forecast algorithm which explains 75\%. However, the increase in stiffness also results in an increase of interpretability. The algorithm allows us to decompose the data into the co- and counter-rotating modes, Figure \ref{fig:008_modes}b and \ref{fig:008_modes}c. This enables practitioners to examine and study modes individually. Figure \ref{fig:008_modes}d shows the prediction of the individual modes, the aggregate prediction (Mode 1 + Mode 2) and ground truth respectively.

These two Koopman models prove to be useful tools in reduced-order modeling for RDE data, giving reconstructions and predictions which provide a number of advantages. Primarily, the models are extremely low-rank, representing the entire wave field in only 2 modes for the case of Figure \ref{fig:008_modes}. The models also prove to be applicable on RDE dynamics that behave linearly and nonlinearly, preserving important nonlinear wave interactions. Another advantage of the Koopman approach is that the methods are robust to noise that is rife in the experimental RDE data, effectively acting as a de-noising filter as can be seen in Figure \ref{fig:387}g. 
Chiefly, the most advantageous aspect of the Spatiotemporal Koopman algorithm is that it is able to separate the wave groups cleanly, and provide a clear representation of the waves traveling in each direction in the case shown in Figure \ref{fig:008_modes}. This is particularly useful in the RDE: the traveling wave shapes and velocities give direct, though qualitative, indication of wave strength, chemical reaction rate, and relative strengths of dissipative effects. Wave strength can be inferred by base-to-peak amplitude of the waves, corresponding to a shock jump condition. Chemical reaction rate can be related to a chemical length scale and is typically observed as the distance from the shock front to the point of greatest luminosity, i.e. the peak of the waveform. Lastly, the rate of decay of the expansion-side of the waves relates to the time scales associated with expelling the burnt combustion products away from the combustion zone. However, there is a lack of interpretability in the wave interactions learned in the Koopman algorithms. This limits the use of Koopman and neural network architecture in gaining intuition into the physics at play in the RDE system. While methods such as UnTWIST may not perform as robustly in determining and separating modes in traveling wave systems, they are able to provide more insight into the underlying dynamics.    

\section{Discussion and Conclusions}

Data-driven ROMs are of growing importance across the engineering, physical and biological sciences given our increasing ability to exploit emerging sensor technologies to observe and quantify complex dynamical systems.  Building models directly from observational data is at the fore-front of data-driven science and engineering~\cite{Brunton2019book} and highlights the 4th-paradigm of data-intensive discovery~\cite{hey2009fourth}.  Importantly, good ROMs require that an appropriate coordinate system be used in order for a low-rank representation of the dynamics to be achieved~\cite{Benner2015siamreview}.  Invariances, particularly rotational and translational, present significant challenges in making ROM models useful for spatiotemporal systems.  Simple traveling waves compromise standard ROM architectures, thus requiring additional methods to handle the translational invariance~\cite{Brunton2019book,kirby1992reconstructing,Rowley2000physd,Rim2018juq,Reiss2018jsc}. More recently, automated methods have been developed to handle traveling waves~\cite{mendible2020dimensionality}.  The so-called UnTWIST method uses spectral clustering and machine learning techniques to provide a reference frame pinned to a traveling wave.  

We have shown that the UnTWIST method can be used on observational data of a rotating detonation engine to find a coordinate system that is pinned to any desired detonation wave.  The robust transformation gives a rotating coordinate system which is ideal for constructing ROMs that characterize the detonation front interactions.  The ROMs are constructed directly from data, requiring no previous physics knowledge of the complex, multi-scale physics driving the combustion dynamics themselves.  Moreover, a diverse of techniques can be applied.  We demonstrated three modeling paradigms: (i) the DMD for building the best-fit linear dynamics model, (ii) a nonlinear Lotka-Volterra model for constructing nonlinear dynamical systems models for the detonation wave interactions, and (iii) a deep Koopman model that uses a neural network to map the time-dynamics to Fourier temporal behavior in order to characterize the dynamics.  All three modeling paradigms are relevant as the RDE data and detonation front interactions exhibit dynamics that range from approximately linear to strongly nonlinear.  Such models provide critical reductions for the complex and multi-scale dynamics of the reactive, compressible fluid dynamics of RDEs.

The overall data-driven architecture emphasis the critical components necessary for physics discovery, specifically the joint discovery of coordinates and parsimonious models that represent interpretable and extrapolatory models of the physics.  Given the recent emergence of RDE data, and the lack of theory characterizing detonation wave interactions, our data-driven method gives both new physical insights, and the first theoretical models of how detonation waves interact.  This allows for engineering design and suggests control strategies that can be imposed in order to manipulate the output of RDE.  This can also inform engineers how to engineering the thermodynamic work loop~\cite{Koch2020b} in order to optimize engine performance.  All of this can be achieved even if a detailed physics model is not available, or if the computations are intractable.

\section*{Acknowledgements}
We acknowledge the support from the Defense Threat Reduction Agency HDTRA1-18-1-0038.  JNK acknowledges support from the Air Force Office of Scientific Research (AFOSR) grant FA9550-17-1-0329.  SLB acknowledges support from the Army Research Office (ARO W911NF-19-1-0045).  JNK and SLB acknowledge support from the National Science Foundation (NSF HDR award \#1934292). 

\begin{spacing}{.8}
 \small{
 \setlength{\bibsep}{4.1pt}
 }

 \end{spacing}
 
\end{document}